\documentclass{article}

\newfont{\bb}{msbm10 at 11pt}
\def\r{\hbox{\bb R}}
\def\l{\hbox{\bb L}}
\def\x{\hbox{\bf x}}

\newenvironment{proof}{\trivlist
\item[\hskip\labelsep{\it Proof}\,:]}{\hfill{$q.e.d.$}\endtrivlist}
\newtheorem{theorem}{Theorem}[section]
\newtheorem{lemma}[theorem]{Lemma}
\newtheorem{proposition}[theorem]{Proposition}
\newtheorem{corollary}[theorem]{Corollary}
\newtheorem{definition}[theorem]{Definition}
\newtheorem{remark}[theorem]{Remark}
\newtheorem{claim}{Claim}

\title{Stationary bands in three-dimensional Minkowski space\thanks{Partially
supported by MEC-FEDER 
 grant no. MTM2004-00109.}} 
\author{Rafael L\'opez\\
Departamento de Geometr\'{\i}a y Topolog\'{\i}a\\
Universidad de Granada\\
18071 Granada (Spain)\\
email:rcamino@ugr.es}
\date{}

\begin{document}

\maketitle
 
\begin{abstract}
In this paper we consider a free  boundary problem in the $3$-dimensional Lorentz-Minkowski space $\l^3$ which deals  spacelike surfaces whose mean curvature is a linear
 function of the time coordinate and the boundary moves in a 
given support plane. We study spacelike surfaces that project one-to-one into a 
strip of the  support  and   that locally  are 
critical points of a certain energy functional involving  the area of the 
surface, a timelike potential and preserves the volume enclosed by the surface. We call 
these surfaces stationary bands. We establish existence of such surfaces and 
we investigate  their  qualitative properties. Finally, we give estimates of its  size  in terms of the initial data. 
\end{abstract}

\section{Introduction and statement of results}\label{intro}

Let $\l^3$ denote the $3$-dimensional Lorentz-Minkowski space, that is, the real vector space $\r^3$ endowed
with the Lorentzian metric $\langle,\rangle=dx_1^2+d_2^2-dx_3^2$, where $x=(x_1,x_2,x_3)$ are the canonical
coordinates in $\r^3$. Let $\Pi$ be a spacelike plane, which we shall assume horizontal,  and 
consider a potential energy $Y$ that, up constants, measures at each point the distance to $\Pi$.   We are interested in the following

\begin{quote}{\it Variational problem.}  To find those spacelike compact surfaces with 
maximal surface area whose boundary moves on 
$\Pi$  and enclosing  a fixed volume of the ambient space. We assume the effect of the potential $Y$.
\end{quote}

The plane $\Pi$ is called the support plane. 
In an approach up to the first order, we are interesting for a spacelike 
compact surface $S$  that is a critical point of  a energy functional  for any 
perturbation of the surface in such way that $S$ is adhered to 
$\Pi$  and the volume determined by $S$ with a 
certain domain of $\Pi$ is prescribed. The energy of 
the system involves the surface area of $S$, the area  of the
domain $\Omega$ in $\Pi$ bounded by the boundary $\partial S$
 of $S$ and the potential defined by $Y$: 
$$E=|S|-\cosh(\beta) |\Omega|+\int_{S} Y\ dS,$$
where $\beta$ is a constant. Consider then an admissible variation of $S$  to our problem, that 
is, a one parameter differentiable family of surfaces $S_t$ indexed by a parameter $t$, 
with $ S_0=S$,  and supported all in $\Pi$: $\partial S_t\subset \Pi$. 
We assume that $S_t$ is a volume preserving variation which means that the value of the 
volume of the region of $\l^3$ determined by $S_t\cup \Omega_t$ is prescribed. We denote by $E(t)$ the corresponding energy of the surface $S_t$. We ask for the shape of those surfaces $S$ that are critical points of the energy for all admissible variations, that is,  
$$\frac{d}{dt}{\Bigg|}_{t=0}E(t)=0.$$
In such case, we  say  that $S$ is a {\it stationary surface}. According to the principle of virtual works, stationary surfaces are characterized by the following:

\begin{theorem}A spacelike surface $S$ in $\l^3$ is stationary 
if  and only if the following two conditions hold: 
\begin{enumerate} 
\item The mean curvature $H$ of  $S$
 is a linear function on  the distance to $\Pi$: 
\begin{equation}\label{statio}
2H(x)=\kappa x_3(x)+\lambda,\hspace*{1cm}\mbox{(Laplace equation)}
\end{equation}
 where $\kappa$ is a constant called the {\it capillary constant}, and $\lambda$ is a constant to be determined by the volume constraint.
\item The surface $S$ intersects $\Pi$ at a constant hyperbolic angle $\beta$  along $\partial S$ (Young condition). 
\end{enumerate}
\end{theorem}

 We refer to \cite{ap,bo,bf,lo1} for more details. In absence of 
the  potential $Y$, the constant $\kappa$ is zero and 
 $S$ is a spacelike surface with constant mean curvature. In this 
sense and in the context of our variational problem, some results have been obtained \cite{alp,ap}. Constant mean curvature spacelike surfaces are well known from the physical point of view because of their role in different problems in General Relativity (see for instance \cite{cy,mt} and references therein). The compact case and $\kappa\not=0$  has been recently considered  by the author in a sucessive 
of works \cite{lo1,lo2,lo3}. The present paper continues 
this work by studying the case that the surface is 
not compact.

In this paper we study spacelike surfaces $S$ that project one-to-one on a strip 
$\Omega_a=\{(x_1,x_2,0);-a<x_1<a\}$ of $\Pi$, 
or in a more general case, in the whole plane $\{x_3=0\}$.  We say then that $S$ is 
a {\it  band} on $\Omega_a$, or simply, a band. Our first motivation comes from the theory of spacelike surfaces with constant mean curvature. In this sense, we point out that maximal bands  ($H=0$) with singularities in $\l^3$ have been studied in the literature \cite{km,mi}. On the other hand, the simplest examples of bands with non-zero
 constant mean curvature are  hyperbolic cylinders: up isometries of the ambient, they are defined by 
${\cal H}_a=\{(\frac{1}{m}\sinh{(x_1)},x_2,\frac{1}{m}\cosh{(x_1)});-a<x_1<a, 
x_2\in\r\}$, $m>0$ and 
the mean curvature is $H=m/2$. A hyperbolic cylinder is also the graph of the 
function $y(x_1,x_2)=\sqrt{x_1^2+1/m^2}$ defined on $\Omega_a$. On the other 
hand, hyperbolic cylinders are surfaces    translationally invariant with respect to a horizontal vector $\vec{w}$. From this point of view, this motivates to study the 
shape of this kind of surfaces, called generalized cylinders and that satisfy the 
Laplace equation (\ref{statio}).  In particular, 
the intersection of $S$ with the given support  plane $\Pi$ realizes with constant 
angle, and so, the Young condition  is satisfied at any boundary point. 

The purpose of this paper is to  study stationary bands establishing existence and certain qualitative features of these surfaces. We begin proving:

\begin{quote} {\it Let $\Omega_a$ be a strip in a 
spacelike plane $\Pi$.  Given $\kappa, \lambda$ and $\beta$ real numbers, there exists 
a stationary band $S$ supported on  $\Omega_a$ that 
satisfies the Laplace equation  and  makes a hyperbolic angle $\beta$ with $\Omega_a$ along its boundary  (Theorems \ref{sessile-ex} and \ref{pendent-ex}). 
Moreover, the surface can extend to be a entire surface.}
\end{quote}

For this, we reduce Equation (\ref{statio}) into an ordinary differential 
equation of second order and we analyze the existence of solutions. This is carried out in 
Section \ref{existence}. The qualitative properties of the shapes that  a 
stationary band adopts depends on  the sign of $\kappa$. We call   {\it sessile} or {\it pendent} stationary 
band if $\kappa>0$ or $\kappa<0$, respectively (this terminology has its origin in the 
Euclidean setting). In Sections \ref{sessile} and 
\ref{estimates} we study the case $\kappa>0$. We prove 

\begin{quote}{\it A sessile stationary band is a convex surface and 
asymptotic to a lightlike cylinder at infinity (Theorem \ref{t-sessile}).}
\end{quote}

Next, we continue studying properties of monotonicity on the parameter $\kappa$ and  we compute the size of the surface in terms of  given data
 in the variational problem.   We omit the 
statements and we refer to Section \ref{estimates} for details. 
Finally in Section  \ref{pendent} 
we study the shape of pendent stationary bands and we describe such surfaces: 

\begin{quote} {\it A pendent stationary  band is invariant by a group of translations 
whose  direction of translation is orthogonal to the rulings of the surface. The time coordinate is a periodic function and the surface extends to an entire spacelike surface 
 (Theorem \ref{t-p}).}
\end{quote}

\section{Preliminaries }\label{preli}

A nonzero vector $v\in \l^3$ is called spacelike or timelike  if 
$\langle v,v\rangle>0$ or $\langle v,v\rangle<0$, respectively. 
Let $S$ be a (connected) surface and let $x:S\rightarrow \l^3$ be an 
immersion of $S$ into $\l^3$. The immersion is said to be spacelike if its tangent 
vectors are spacelike. Then the scalar product $\langle,\rangle$ induces 
a Riemannian metric on $S$. 
Observe that $\vec{e_3}=(0,0,1)$ is a unit timelike vector field globally defined on $\l^3$, which determines 
a time-orientation on the space $\l^3$. This allows us to choose a unique unit normal vector field $N$ on $S$ which is in the 
same time-orientation as $\vec{e_3}$, and hence that $S$ is oriented by $N$. In this 
article  all spacelike surfaces  will be oriented according to this choice of $N$. Because 
the support plane in our variational problem is horizontal,  the hyperbolic angle $\beta$ between $S$ and $\Pi$ along its boundary is given by
$\langle N,\vec{e_3}\rangle=-\cosh\beta$.

For spacelike immersions, the notions of 
the first and  second fundamental form are defined in the same way as  in Euclidean space, 
namely, 
$${\rm I}=\sum_{i j}g_{ij} dx_i\ dx_j,\hspace*{.5cm}\mbox{and}\hspace*{.5cm}{\rm II}=\sum_{i j} h_{ij} dx_i\ dx_j,$$
respectively, where $g_{ij}=\langle\partial_i x,\partial_j x\rangle$ is the induced metric on $S$ by $x$ and 
$h_{ij}=\langle \partial_i N,\partial_j x\rangle$.
Then the  mean curvature $H$ of $x$ is given  by 
\begin{equation}\label{mean}
2H=\mbox{trace\ ($I^{-1} II$)}=\frac{h_{11} g_{22}-2h_{12} 
g_{12}+h_{22}g_{11}}{\mbox{det}(g_{ij})}.
\end{equation}
Locally, if we  write  $S$ as  the graph of a smooth function $u=u(x_1,x_2)$ 
defined over a domain $\Omega$, the spacelike condition implies $|\nabla u |<1$. 
According to the choice of the time orientation, $N$ is 
$$N=\frac{(\nabla u,1)}{\sqrt{1-|\nabla u|^2}},$$
 and the mean curvature $H$ of $S$ at each point 
$(x,u(x))$,  $x\in\Omega$,  satisfies the equation
$$(1-|\nabla u|^2)\Delta u+\sum u_i u_j u_{ij}=2H(1-|\nabla u|^2)^{3/2}.$$
This equation is of quasilinear elliptic type and it 
can alternatively be written in divergence form
\begin{equation}\label{media}
\mbox{div}\left(\frac{\nabla u }{\sqrt{1-|\nabla u |^2}}\right)=2 H.
\end{equation}
In particular, if $u$ and $v$ are two  functions are solutions of the same 
equation (\ref{media}),  the difference function $w=u-v$ satisfies an 
elliptic linear equation $Lw=0$ and  one can apply the Hopf maximum principle. 
Then we obtain uniqueness of solutions for each given boundary data.

We now consider the type of surfaces which  are interesting in this work, and that
 generalize the family of hyperbolic cylinders defined in the Introduction.
 These surfaces are cylindrical surfaces, also called in the literature, generalized cylinders. A cylindrical surface 
$S$ is a ruled surface generated by a one-parameter family of straight-lines 
$\{\alpha(s)+t\vec{w}; t\in \r\}$, parametrized by the parameter $s$,  where $\alpha(s)$, $s\in I$, is a regular curve contained 
in a plane $P$ and  $\vec{w}$ is a given vector which is not parallel to $P$. The 
curve $\alpha$ is called a directrix of $S$ and the lines are called the rulings. 
The shape of a cylindrical surface is completely determined then by the 
geometry of $\alpha$. In addition, if we impose that $S$ is a spacelike surface, 
 then  both $\alpha'(s)$ and 
$\vec{w}$ are spacelike vectors. 
For example, a  hyperbolic cylinder is a cylindrical surface in 
$\l^3$ whose directrix is, up isometries of the ambient space,  a (spacelike) hyperbola 
in a vertical plane $P$ and $\vec{w}$ is a horizontal vector orthogonal to $P$. 

In a first moment, 
the directrix $\alpha$ 
can have self-intersections and so, the surface determined 
could not be a band on $\Pi$.  However, if we impose  
the condition to be stationary, then 
 we prove that the 
rulings of the surface must be horizontal, and then, $\alpha$ is an embedded curve. Exactly:

\begin{proposition}\label{32} Let $S$ be a spacelike cylindrical  surface in $\l^3$. If the 
mean curvature of $S$ is a linear function of the time coordinate with 
$\kappa\not=0$, then the 
rulings  are horizontal.
\end{proposition}

\begin{proof} We parametrize $S$  as 
$$\x(s,t)=\alpha(s)+t \vec{w},\hspace*{1cm}s\in I, t\in\r,$$
where $\alpha$ is the directrix of $S$  parametrized by the length arc. 
The computation of the mean curvature $H$ of $S$ according to (\ref{mean}) gives
\begin{equation}\label{2h}
2H=\frac{\langle\alpha'(s)\times\vec{w}, 
\alpha''(s)\rangle}{(1-\langle\alpha'(s),\vec{w}\rangle^2)^{3/2}},
\end{equation}
where $\times$ indicates the cross product in $\l^3$. In particular, the mean 
curvature function depends only on the $s$-variable. 
Therefore, if $H$ 
satisfies the relation (\ref{statio}), that is, 
$$H=\kappa x_3\circ(\alpha(s)+t\vec{w})+\lambda=\kappa x_3\circ\alpha(s)+\lambda+
\kappa (x_3\circ\vec{w})\ t,$$
 we infer that 
$x_3\circ \vec{w}=0$, and then, $\vec{w}$ is a horizontal vector.
\end{proof}

As a consequence of Proposition \ref{32}, we can choose the plane 
$P$ containing the directrix  to be vertical,  with the rulings being 
horizontal straight-lines.  On the 
other hand, and because $\alpha$ is a spacelike planar curve in a vertical plane, 
 $\alpha$ is globally the graph of a certain function $u$ defined 
on an interval of any horizontal line of $P$. Hence, we conclude:

\begin{corollary} Any   cylindrical surface in $\l^3$ 
 that satisfies the Laplace equation (\ref{statio}) is a  band. 
\end{corollary}

\begin{definition} A stationary band in $\l^3$ is a cylindrical surface that 
satisfices the Laplace equation (\ref{statio}).
\end{definition}

\section{Stationary bands: existence and symmetries}\label{existence}

In this section, we write our variational problem in terms of the theory of 
ordinary differential equations, exactly, the Laplace equation (\ref{statio}) reduces in an  
ordinary differential equation  of second order. 
The purpose of this section is to establish results of existence of the corresponding boundary value problem together with properties of 
symmetries of the solutions. 
Let $S$ be a stationary band on a strip $\Omega_a$  of $\l^3$. We have then
$$\kappa\ x_3+\lambda=\langle\alpha'\times\vec{w},\alpha''\rangle=C_{\alpha},$$
where $C_{\alpha}$ denotes the curvature of $\alpha$.  In addition, 
the angle that makes a such surface with the given horizontal support 
plane  is constant. 

If $\kappa=0$ in the Laplace equation (\ref{statio}), the surface has  constant mean curvature  $H=\lambda/2$. Then the curvature $C_{\alpha}$ is constant, namely, $C_{\alpha}=\lambda$. Therefore $\alpha$ is a straight-line or  a 
spacelike hyperbola of $\l^3$ and the corresponding  
surfaces are planes or hyperbolic cylinders respectively. Assuming that $\kappa\not=0$,
 we do a change of variables to get $\lambda=0$ in the Laplace equation. For this, 
it suffices with the change of 
the immersion $x$ by $x:=x+(\lambda/\kappa)\vec{e_3}$. Then the new 
spacelike surface is a stationary band on the same strip $\Omega_a$ and with 
 mean curvature 
$H(x)=\kappa x_3$. The support plane $\Pi$ has changed by other horizontal plane, namely, $\{x_3=\lambda/\kappa\}$ and that we will denote by the same symbol $\Pi$. However, the 
contact agle of $x$  along its boundary is 
$\beta$ again. It follows that the shape of the original surface 
$S$ is independent of the constraint $\lambda$ in (\ref{statio}). 
Throughout  this work, 
we shall consider that $\kappa\not=0$ and that $\lambda=0$ in the Laplace equation 
(\ref{statio}). Then the problem of existence of our variational problem re-states as follows: 

\begin{quote} {\it Variational problem.} Let 
$\Omega_a $ be a strip of the  $(x_1,x_2)$-plane and let $\beta$ and 
$\kappa\not=0$ two real numbers. Does there exists  a stationary band $S$ on $\Omega_a$ such that: i) its mean curvature is $H(x)=\kappa x_3$ in each point $x\in S$ and; ii) the angle that makes 
$S$ with a horizontal plane $\Pi$ along $\partial S$ is $\beta$?
\end{quote}

Let us consider $S$ a stationary band defined on the strip $\Omega_a$ and we parametrize 
by $S=\{r,x_2,u(r));-a<r<a,x_2\in \r\}$. According to the choice of the orientation  on $S$, 
 the hyperbolic angle $\beta$ between $S$ and $\Pi$ along $\partial S$ satisfies
$$\cosh\beta=-\langle N_S,\vec{e_3}\rangle=-\langle(\frac{(u_{x_1},u_{x_2},1)}{\sqrt{1-u'^2}},
(0,0,1)\rangle=
\frac{1}{\sqrt{1-u'^2}}\hspace*{1cm}\mbox{at $|r|=a$}.$$
In the variational problem,  the Laplace  and 
the Young equations  write, respectively, as 
\begin{eqnarray}
 & & \frac{ u''(r)}{(1-u'(r)^2)^{3/2}}=\kappa\ u(r),\hspace*{1cm}-a< r<a.\label{f1}\\
& & u'(\pm a)=\pm\tanh\beta.\label{f2}
\end{eqnarray}

We change (\ref{f1})-(\ref{f2})
 by the initial value problem
\begin{eqnarray}
& &\frac{ u''(r)}{(1-u'(r)^2)^{3/2}}=\kappa\ u(r),\hspace*{1cm} r>0\label{e1}\\
& &u(0)=u_0,\hspace*{1cm}u'(0)=0\label{e2}
\end{eqnarray}
where $u_0$ is a real number.

\begin{theorem} Given $u_0$, there exists a unique solution of (\ref{e1})-(\ref{e2}).
The solution $u=u(r;u_0,\kappa)$ depends continuosly  on the parameters $u_0$, $\kappa$ and the maximal interval of definition of $u$ is $\r$.
\end{theorem}

\begin{proof}
Put $v=u'/\sqrt{1-u'^2}$. 
Then the problem (\ref{e1})-(\ref{e2}) becomes to 
equivalent to a pair of differential equations
\begin{eqnarray}
u'&=&\frac{v}{\sqrt{1+v^2}},\hspace*{.5cm}u(0)=u_0.\label{s1}\\
v'&=&\kappa u,\hspace*{1.5cm}v(0)=0.\label{s2}
\end{eqnarray}
The solution $u$ that we look for is then defined by 
\begin{equation}\label{uve}
u(r)=u_0+\int_0^r \frac{v(t)}{\sqrt{1+v(t)^2}}\ dt.
\end{equation}
Then standard existence theorems of ordinary differential equations assures local existence and uniqueness of (\ref{s1})-(\ref{s2}) and the continuity of 
solutions with respect to the parameters $u_0$ and $\kappa$. 
We study the maximal domain of the solution.  
On the contrary, suppose that $[0,R)$ is the maximal interval 
of the solution $u$, with $R<\infty$. 
By (\ref{s2}) and (\ref{uve}), 
$$|u(r)|<|u_0|+r,\hspace*{1cm}|v(r)|\leq \kappa r(|u_0|+\frac{r}{2}).$$
 Thus, and using (\ref{s1})-(\ref{s2}), the limits of $u'$ and $v'$ at $r=R$ are finite, which it would imply that we can extend the solutions $(u,v)$ beyond $r=R$: contradiction.
\end{proof}

\begin{corollary} Any stationary band of $\l^3$ supported on a spacelike plane $\Pi$ 
 can extend to be a graph defined in the whole plane, that is, it is an entire 
spacelike surface of $\l^3$. 
\end{corollary}

If the constant $\kappa$ is understood we write  by $u(r;u_0)$ a 
 solution of the initial value problem (\ref{e1})-(\ref{e2}).

It is possible to obtain a 
first integration of (\ref{e1})-(\ref{e2}) multiplying by $u'$ in (\ref{e1}): 
\begin{equation}\label{cotaz}
u^2=u_0^2+\frac{2}{\kappa}\left(\frac{1}{\sqrt{1-u'^2}}-1\right).
\end{equation}

Therefore
$$u(r)=u_0+\int_0^r\sqrt{1-\frac{4}{(2+\kappa(u(t)^2-u_0^2)^2}}\ dt.$$
Denote $\psi=\psi(r)$ the hyperbolic angle that makes the 
directrix $\alpha(r)=(r,u(r))$ with the horizontal direction. Put
\begin{equation}\label{angulo}
\sinh\psi=\frac{u'}{\sqrt{1-u'^2}},\hspace*{1cm}\cosh\psi=\frac{1}{\sqrt{1-u'^2}},
\end{equation}
Then the Euler-Lagrange equation (\ref{e1}) takes the form 
\begin{equation}\label{B1}
(\sinh\psi)'=\kappa u, 
\end{equation}
and so, 
\begin{equation}\label{B11} 
\sinh\psi=\kappa \int_0^r u(t)\ dt.
\end{equation}
The identity (\ref{B1}) corresponds actually with the mean curvature equation  
in its divergence form (\ref{media}). Using (\ref{cotaz}), Equation (\ref{statio}) writes now as 
\begin{equation}\label{cotaz2}
u^2=u_0^2+\frac{2}{\kappa}\left(\cosh\psi-1\right).
\end{equation}
We  study the symmetries of the shape of a stationary band. Of course, each  plane orthogonal to the rulings is a plane of symmetry of $S$. However, stationary bands in $\l^3$ have a rich symmetry. 

\begin{theorem}[Symmetry]\label{symmetry} Let $u$ be a solution of (\ref{e1}). 
\begin{enumerate}
\item  If $u'(r_0)=0$,  the graphic of $u$ is symmetric with respect to 
 the vertical line $\{r=r_0\}$.
\item If  $u(r_0)=0$,  the graphic of $u$ is symmetric with respect 
to the point $(r_0,0)$. 
\end{enumerate}
\end{theorem}
 
\begin{proof} In both cases, we can assume that $r_0=0$. We prove the first statement, that is, that $u(r)=u(-r)$. 
The functions $u(r)$ and $u(-r)$ are solutions of the same equation (\ref{e1}) and with the 
same initial conditions at $r=0$, namely, $u_0$ and $u'(0)=0$. Then the uniqueness of 
solutions gives $u(r)=u(-r)$.   
The proof of the second statement is similar in showing $u(r)=-u(-r)$: now both functions are solutions of (\ref{e1}) with initial conditions $u_0=0$ and $u'(0)$.
\end{proof}

Finally, we establish a result that says us that the sign of $u_0$ and 
$\kappa$ can take the same one.

\begin{proposition} Let $u=u(r;u_0,\kappa)$ be a solution of (\ref{e1})-(\ref{e2}). 
Then $u(r;u_0,\kappa)=-u(r;-u_0,\kappa)$.
\end{proposition}

\begin{proof} Again, this is a direct consequence of the uniqueness of solutions.
\end{proof}

As a consequence of this result, one can choose $u_0$ to have the same sign than $\kappa$. This will be assumed 
 throughout the text.

\section{Stationary  bands: the case $\kappa>0$}\label{sessile}

This section is devoted to study the qualitative properties of the shape of a 
sessile stationary band. Assume $\kappa>0$. Recall that  $u_0>0$.
The geometry of the directrix of a sessile stationary band is described by the 
next:

\begin{theorem}[Sessile case]\label{t-sessile} Let $u=u(r;u_0)$ be a solution of 
(\ref{e1})-(\ref{e2}).  
Then the function $u$ has exactly a minimum at $r=0$  with 
$$\lim_{r\rightarrow\infty}u(r)=\infty,\hspace*{1cm}\lim_{r\rightarrow\infty}u'(r)=1.$$
Moreover, $u$ is convex with
$$\lim_{r\rightarrow\infty}u''(r)=0.$$
\end{theorem}

\begin{proof} 
Since $u'(0)=0$,  we restrict the study $u$ for $r\geq 0$ (Theorem \ref{symmetry}). 
The integrand in (\ref{B11}) is positive near to $r=0$. Then $\sin\psi>0$, and so, $u'(r)>0$.
This means that $u$ is increasing near $r=0$. If $r_o$ is the first point where 
$u'(r_o)=0$, then (\ref{cotaz}) implies that $u(r_o)=u(0)$: contradiction. Thus, 
$u'(r)>0$ for any $r$ and this proves that $u$ is strictly increasing 
and $r=0$ is the only minimum. On the 
other hand, at $r=0$, $u''(0)=\kappa u_0>0$, which implies that $u$ is convex around $r=0$. Since $u(r)>u_0>0$, Equation (\ref{e1}) 
concludes that $u''$ has not zeroes, that is, $u$ is a convex function. 

As $u(r)\rightarrow\infty$ as $r\rightarrow\infty$, it follows from (\ref{cotaz}) that 
$\cosh\psi(r)\rightarrow\infty,$
that is, $u'(r)\rightarrow 1$ as $r\rightarrow\infty$. Finally, from (\ref{e1}) and 
(\ref{cotaz2}), 
$$0\leq u''=\frac{\kappa}{\cosh^3\psi}\sqrt{u_0^2+(\cosh\psi-1)}
\rightarrow 0,$$
as $r\rightarrow\infty$.
\end{proof}

As consequence of Theorem \ref{t-sessile},
 any sessile stationary band has a lowest  height if $u_0>0$ (or 
upper height if $u_0<0$) and this height is reached at the ruling 
 $\{(0,x_2,u_0);x_2\in\r\}$. 
The fact that $u'(r)\rightarrow 1$ as $r\rightarrow\infty$ can be written as follows: 

\begin{corollary} Any sessile stationary band of  $\l^3$ is asymptotic at infinity to a 
lightlike cylinder of $\l^3$.
\end{corollary}

Recall that, up isommetries, a lightlike cylinder is the surface of $\l^3$
 defined  as $\{(x_1,x_2,x_3)\in\l^3;x_1^2-x_3^2=0\}$.

As $u$ is increasing on $(0,\infty)$, we bound the integrand in (\ref{B11}) by $u_0<u(t)<u(r)$ obtaining 
\begin{equation}\label{B2}
 \kappa   u_0<\frac{\sinh\psi(r)}{r}<\kappa  u(r).
\end{equation}
Moreover, 
\begin{equation}\label{u0}
\lim_{r\rightarrow 0}\frac{\sinh\psi(r)}{r}=\kappa u_0.
\end{equation}

We establish the existence  of the original variational problem introduced 
in this work.  

\begin{theorem}[Existence]\label{sessile-ex}  Let $\Omega_a$ be a strip 
of the $(x_1,x_2)$-plane, $a>0$.  
Given constants $\kappa>0$ and $\beta$,  there exists a stationary band on 
$\Omega_a$ whose directrix is defined 
by a function $u=u(r;u_0)$, that makes a contact hyperbolic angle 
$\beta$ with the support plane $\{x_3=u(a)\}$. 
\end{theorem}

\begin{proof} 
If $\beta=0$, we take $S=\{x_3=0\}$. Without loss of generality, we now assume 
$\beta>0$. The problem is equivalent to search a solution of (\ref{f1})-(\ref{f2}). For this, we take the initial value problem (\ref{e1})-(\ref{e2}), with $u_0>0$. The 
problem then reduces to find $u_0>0$ such that $u'(a;u_0)=\tanh\beta$. By the continuity of 
the paramaters, $\lim_{u_0\rightarrow 0}u'(a;u_0)=0$. On the other hand, 
by using (\ref{B2}), 
$$u'(a;u_0)=\tanh\psi(a)\geq\frac{\kappa u_0 a}{\sqrt{1+\kappa^2 u_0^2 a^2}}\longrightarrow 1$$
as $u_0\rightarrow\infty$. 
As a consequence, and by continuity again, given $\beta>0$, we can find $u_0>0$ such that 
$u'(a;u_0)=\tanh\beta$.
\end{proof}

The uniqueness of solutions will be derived at the end of this section: 
see Corollary \ref{uni}. 
After the existence of solution of (\ref{f1})-(\ref{f2}), we continue showing certain results regarding the dependence of the solutions on the parameters of the differential equation, as 
for example, $\kappa$ and $u_0$. Exactly, we establish results about the monotonicity. 
First, we concentrate on the capillary constant $\kappa$ and we begin by proving the next result about the 
solutions of (\ref{e1})-(\ref{e2}).

\begin{theorem}\label{mono-k} 
Let $\kappa_1,\kappa_2>0$. Denote $u_i=u_i(r;u_0,\kappa_i)$, $i=1,2$,  two solutions of (\ref{e1})-(\ref{e2}) with the 
same initial condition $u_0$.   If $\kappa_1<\kappa_2$, then 
$u_1(r)<u_2(r)$ for any $r\not=0$ and $u_1'(r)<u_2'(r)$ for $r>0$.
\end{theorem}

\begin{proof} Denote by $\psi^{(i)}$ the angle functions defined by 
(\ref{angulo})  for each function $u_i$. 
We know from (\ref{B11}) that
\begin{equation}\label{p-mono}
\sinh\psi^{(2)}(r)-\sinh\psi^{(1)}(r)=\int_0^{r}(\kappa_2 u_2(t)-\kappa_1 u_1(t))\ dt.
\end{equation}
At $r=0$, the integrand is positive and so, $\psi^{(2)}(r)>\psi^{(1)}(r)$ on 
some interval $(0,\epsilon)$. Then $u_2'(r)>u_1'(r)$ and because 
$u_2(0)=u_1(0)$, we have $u_2(r)>u_1(r)$ in $(0,\epsilon)$. We prove that $u_2'(r)>u_1'(r)$ 
holds for any $r>0$.  If $r_0>0$ is the first point where $u_2'(r_0)=u_1'(r_0)$, then 
$u_2''(r_0)\leq u_1''(r_0)$ and $u_2(r_0)>u_1(r_0)$. But (\ref{e1}) gives
$$u_2''(r_0)=C_{u_2}(r_0)=\kappa_2 u_2(r_0)>\kappa_1 u_1(r_0)=C_{u_1}(r_0)=u_1''(r_0).$$
This contradiction implies that $u_2'(r)>u_1'(r)$ for any $r>0$ and then, $u_2(r)>u_1(r)$.

\end{proof}

We now return to the boundary value problem (\ref{f1})-(\ref{f2}). We prove that
 for fixed $\beta$ and $a$, the solution 
of (\ref{f1})-(\ref{f2}) and its derivative with respect to $r$ are monotone functions 
of $\kappa$. 

\begin{theorem}[Monotonicity with respect to $\kappa$] \label{mono} Let $\kappa_1, 
\kappa_2>0$. Denote $u_i=u_i(r)$, $i=1,2$,  two solutions of (\ref{f1})-(\ref{f2}) on 
the strip $\Omega_a$ for $\kappa=\kappa_i$ with 
$u_i(0)>0$. If $\kappa_1<\kappa_2$, then 
\begin{enumerate}
\item $u_1(r)>u_2(r)$ for $0\leq r\leq a$.
\item  $u_1'(r)>u_2'(r)$ for $0<r<a$. 
\end{enumerate}
\end{theorem}

\begin{proof} Put $v_i=\sinh\psi^{(i)}$.  For 
each $0\leq r_0<r<a$, (\ref{B11}) yields:
$$v_i(r)- v_i(r_0)=\int_{r_0}^r  \kappa_i u_i(t)\ dt.$$
Then  
\begin{equation}\label{mono1}
v_2(r)-v_1(r)=v_2(r_0)-v_1(r_0)+\int_{r_0}^r (\kappa_2 u_2(t)-\kappa_1 u_1(t)) \ dt.
\end{equation}
Because $u_1'(a)=u_2'(a)$,
\begin{equation}\label{mono2}
v_1(r_0)-v_2(r_0)=\int_{r_0}^a (\kappa_2 u_2(t)-\kappa_1 u_1(t))\ dt.
\end{equation}

\begin{claim}  If   $\kappa_2 u_2(r_0)\geq \kappa_1 u_1(r_0)$, then 
$v_2(r_0)<v_1(r_0)$.
\end{claim}

On the contrary case, that is, if $v_1(r_0)\leq v_2(r_0)$, then $0<u_1'(r_0)<u_2'(r_0)$. 
The fact that  $\kappa_1<\kappa_2$ implies
 $\kappa_1 u_1'(r_0)<\kappa_2 u_2'(r_0)$. Hence $\kappa_1 u_1<\kappa_2 u_2$ on a certain 
interval $(r_0,r_0+\delta)$. Let $r_1\in (r_0,a]$ be the largest number where such inequality holds. 
In view of (\ref{mono1}), 
for each $r_0<s\leq r_1$,  $v_2(s)>v_1(s)$. Thus $u_2'(s)>u_1'(s)$ and 
$\kappa_2 u_2'> \kappa_1 u_1'$. This implies $\kappa_2 u_2>\kappa_1 u_1$ for each $r_0<s\leq r_1$. 
Since $r_1$ is maximal, then $r_1=a$. 
We put now $s=a$ in (\ref{mono2})  and we obtain
$$0\geq v_1(r_0)-v_2(r_0)=\int_{r_0}^a (\kappa_2 u_2(t)-\kappa_1 u_1(t))\ dt>0,$$
 which it is a contradiction. This proves the Claim.

Let us prove now the Theorem and we begin with the item 2. Assume there exists $r_0$, $0<r_0<a$, such that $u_1'(r_0)\leq u_2'(r_0)$. 
Then $v_1(r_0)\leq v_2(r_0)$. By the Claim,  $\kappa_2 u_2(r_0)<\kappa_1 u_1(r_0)$, 
and so, (\ref{e1}) implies $v_2'(r_0)<v_1'(r_0)$. For a 
certain neighbourhood on the left of $r_0$, we obtain then
$$0\leq v_2(r_0)-v_1(r_0)<v_2(r)-v_1(r)$$
 which it yields $v_2(r)>v_1(r)$. As 
$v_1(0)=v_2(0)=0$, there exists a last number $r_1$, $0\leq r_1<r_0$, such that 
$v_2>v_1$ in the interval $(r_1,r_0)$ and $v_2(r_1)=v_1(r_1)=0$. The Claim implies now 
$\kappa_2 u_2(r)<\kappa_1 u_1(r)$, for $r_1<r\leq r_0$. But (\ref{mono2}) yields 
$v_2(r)<v_1(r)$ and that is a contradiction. Consequently,  $u_2'<u_1'$ in $(0,a)$.

Let us prove the item 1. As $u_2'<u_1'$, $v_2<v_1$. As close  $r=0$, 
$\kappa_i u_i(r) r\geq \kappa_i u_i(0) r$, it follows from (\ref{u0}) that
$$\lim_{r\rightarrow 0} \frac{v_i(r)}{r}=\kappa_i u_i(0).$$
Because $v_1>v_2$, we infer then $u_2(0)<u_1(0)$. As $u_2'<u_1'$, an integration leads to $u_2<u_1$ on the interval $[0,a]$.
\end{proof}

 Given a capillary 
constant $\kappa$, we would like to control the dependence
 of solutions $u(r;u_0)$ with respect to the 
initial condition $u_0$. We will obtain monotonicity, that is, if $u_0<v_0$, then $u(r;u_0)<u(r;v_0)$ for any $r$. Moreover, we 
can precise the distance between the two solutions.

\begin{theorem}\label{delta} 
Fix   $\kappa>0$. If $\delta>0$, then $u(r;u_0+\delta)-\delta>u(r;u_0)$ for any $r\not=0$.
\end{theorem}

\begin{proof} 
By symmetry, it is suffices to  show the inequality for $r>0$.
Define the function $u_{\delta}=u(r;u_0+\delta)$, and let $\psi^{\delta}$ be it the corresponding hyperbolic angle, see 
(\ref{angulo}). It follows from (\ref{B11}) that 
\begin{equation}\label{integrand}
\sinh\psi^{\delta}-\sinh\psi=\kappa \int_0^r  (u_{\delta}(t)-u(t))\ dt.
\end{equation}
Since the integrand is positive at $r=0$, there exists $\epsilon>0$ such that 
$$\sinh\psi^{\delta}(r)-\sinh\psi(r)>0 \hspace*{1cm}\mbox{in }(0,\epsilon).$$
As
$$\sinh\psi^{\delta}(0)-\sinh\psi(0)=0,$$
we have the inequality $\psi^{\delta}>\psi$ in the interval $(0,\epsilon)$. In addition, 
$$(u_{\delta}(r)-u(r))'=\tanh\psi^{\delta}(r)-\tanh\psi(r)>0.$$
Therefore the function $u_{\delta}-u$ is strictly increasing on $r$. So,  
$u_{\delta}(r)-\delta>u(r)$. Let $r_0>\epsilon$ be the first point where $u_{\delta}(r_0)-\delta=u(r_0)$. Again 
(\ref{integrand}) yields $\sinh\psi^{\delta}(r_0)-\sinh\psi(r_0)>0$ and $(u_{\delta}-u)'(r_0)\leq 0$.  But this implies that $\psi^{\delta}(r_0)\leq \psi(r_0)$, which  is a contradiction. 
As conclusion, $u_{\delta}-\delta>u$ in $(0,\infty)$ and this shows the result.
\end{proof}

\begin{corollary} Let $S_1$, $S_2$ be two sessile stationary bands on $\Omega_a$ with the same 
capillary constant $\kappa$. Let $h_i$ be the lowest heights of $S_i$, $i=1,2$. 
If $0<h_1<h_2$, we can move $S_2$ by translations until  touches $S_1$ and in such way that  
$S_2$ lies 
completely above $S_1$.
\end{corollary}

\begin{corollary}[Uniqueness] \label{uni} 
The solution obtained in Theorem \ref{sessile-ex} is unique.
\end{corollary}

\begin{proof} By contradiction, assume that $S_1$ and $S_2$ are two 
different stationary bands on $\Omega_a$ and with the 
same Young condition. By the symmetries of solution of (\ref{f1}) and
 Theorem \ref{t-sessile}, $S_1$ and $S_2$ 
are determined by functions $u_1=u(r;u_0) $ and $u_2=u(r;v_0)$ respectively,
 solutions of (\ref{e1})-(\ref{e2}) 
on the same strip $\Omega_a$ and with $u_1'(a)=u_2'(a)$. Without loss of generality, 
we assume 
that $0<u_0<v_0$. The proof of Theorem \ref{delta} says that $u_2'(r)>u_1'(r)$ in some interval $(0,\epsilon)$. Actually, we now prove  that this inequality holds for any $r>0$. If $r_0$ is the first point 
where $u_1'(r_0)=u_2'(r_0)$, then $u_2''(r_0)\leq u_1''(r_0)$ and 
$u_2>u_1$ on $[0,r_0]$. But Eq. (\ref{e1}) implies that 
$$u_2''(r_0)=\kappa u_2(r_0)>\kappa u_1(r_0)=u_1''(r_0).$$
This contradiction implies that $u_2'>u_1'$ for any $r>0$. But then it is impossible 
that $u_2'(a)=u_1'(a)$.
\end{proof}

This section ends with a result of foliations of the ambient space $\l^3$ by   stationary bands.

\begin{corollary} \label{foliation}Fix $\kappa>0$. Then the Lorentz-Minkowski space $\l^3$ can be foliated by 
a one-parameter family of sessile stationary bands, 
for the same capillary constant $\kappa$.  The foliations are given by stationary 
bands that are entire spacelike surfaces whose 
profile curves are $\{u(r;u_0);u_0\in\r\}$ and $u_0$ is the parameter of the foliation. 
\end{corollary}

\begin{proof}
Let $(a,b)$ a point in the $(x_1,x_2)$-plane. We have to show that there exists a unique $u_0$ such that 
$b=u(a;u_0)$. If $b=0$, we take $u_0=0$. We assume  that $b>0$ (the reasoning is similar if $b<0$). By the 
dependence of solutions of (\ref{e1})-(\ref{e2}) and because $u(r;u_0)\geq u_0$ (assuming $u_0>0$), we have
$$\lim_{u_0\rightarrow 0}u(a;u_0)=0,\hspace*{1cm}\lim_{u_0\rightarrow\infty}u(a;u_0)=\infty.$$
Then we employ again the dependence of the parameter $u_0$ to assure the existence of $u_0$ such $u(a;u_0)=b$. The uniqueness of $u_0$ is given by Theorem \ref{delta}.
\end{proof}

\begin{remark} It is worth to point out that in relativity theory 
there is interest of finding 
 real-valued functions on a given spacetime, all of whose level sets provide a global time coordinate. Consequently, 
Corollary \ref{foliation} says that we can find a foliation of
 the ambient space $\l^3$ whose leaves are   sessile stationary bands tending to 
lighlike cylinders at infinity. 
\end{remark}

\begin{remark} In the proofs of the about results, as for example, 
 Theorems \ref{mono-k}, \ref{mono} and \ref{delta},
we have needed to compare solutions of 
equations of type (\ref{e1}). Actually,   
we are using the Hopf maximum principle for 
the mean curvature equation (\ref{media}), but the fact that this equation reduces into one variable, it makes easier the comparisons that we we look for.

\end{remark}

\section{Sessile stationary bands: estimates}\label{estimates}

This section is devoted to obtain estimates of the size for a  sessile stationary band. Exactly, we will give bounds of the height of the solutions of our variational problem in 
terms of the lowest height $u_0$,  the hyperbolic angle of contact 
$\beta$ or the width $2a$ of the strip $\Omega_a$.

Fix $\kappa>0$. Let $S$ be a stationary band on a strip $\Omega_a$ given by a solution 
$u=u(r;u_0)$ of (\ref{e1})-(\ref{e2}) 
and such that $u'(a)=\tanh\beta$. A first control of $u(a)$ comes from the 
identity (\ref{cotaz}): 
\begin{equation}\label{cotaz22}
u(a)=\sqrt{u_0^2+\frac{2}{\kappa}(\cosh\beta-1)}.
\end{equation}
The estimates that we shall obtain are a consequence of the comparison
 of our stationary bands   with hyperbolic cylinders. 
 Consider $y_1$ and $y_2$ two  hyperbolas defined on $[0,a]$ and given by 
\begin{equation}\label{y1}
y_1(r)=u_0-\mu_1+\sqrt{r^2+\mu_1^2},\hspace{1cm}\mu_1=\frac{1}{\kappa u_0}.
\end{equation}
\begin{equation}\label{u2}
y_2(r)=u_0-\mu_2+\sqrt{r^2+\mu_2^2},\hspace*{1cm}
\mu_2= \frac{a}{\sin\psi(a)}.
\end{equation}
Let $\Sigma_1$ and $\Sigma_2$ be the  hyperbolic cylinders obtained by 
translating $y_1$ and 
$y_2$ in the $x_2$-direction, respectively. Both surfaces have constant mean curvature:
$$H_1=\frac{C_{y_1}}{2}=\frac{\kappa u_0}{2},\hspace*{1cm}H_2=\frac{C_{y_2}}{2}= 
\frac{\sinh\psi(a)}{2a}.$$ 
The functions $y_1$ and $y_2$ agree with $u$ at $r=0$. On the other hand, 
the mean curvature of $S$ at the point $(0,x_2,u_0)$ agrees with the one of $\Sigma_1$ and 
$u'(a)=y_2'(a)$.

\begin{lemma} \label{lemma-claim} 
The surface $S$ lies between  $\Sigma_1$ and  $\Sigma_2$.
\end{lemma}

\begin{proof} In order to prove the result, it suffices to show $y_1<u<y_2$  on the interval $(0,a]$. Denote $C_u$ the curvature of the graphic of $u$. 
At $r=0$, 
$C_u(0)=\kappa u_0=C_{y_1}(0)$, but $C_u$ is increasing on $r$ since both $\kappa$ and 
$u'$ are positive.  Because
$u(0)=y_1(0)$, we conclude that $y_1(r)<u(r)$ in 
$0<r<a$. We now prove the inequality $u<y_2$. As the curve $y_2$ 
has constant curvature, inequalities  (\ref{B2}) yield
$$C_{y_2}(0)=C_{y_2}(a)=\frac{\sinh\psi(a)}{a}> \kappa u_0=C_{u}(0).$$
Then at $r=0$, $C_{u}(0)<C_{y_2}(0)$. As $y_2(0)=u(0)$, it follows that 
 there exists $\delta>0$ such that $u(r)<y_2(r)$ for $0<r<\delta$. We assume that 
$\delta$ is the least upper bound of such values. By contradiction, suppose that $\delta<a$. As $y_2(\delta)=u(\delta)$
and $y_2'(\delta)\leq u'(\delta)$, $\psi^{(2)}(\delta)\leq\psi(\delta)$ and thus
\begin{equation}\label{alfa}
\int_0^{\delta}\frac{d}{dr}\left(\sinh\psi(r)-\sinh\psi^{(2)}(r)\right)dr=\sinh\psi(\delta)-\sinh\psi^{(2)}(\delta):=\alpha(\delta)\geq 0.
\end{equation}
Hence there exists $\bar{r}\in (0,\delta)$ such that 
$$C_u(\bar{r})=(\sinh\psi)'(\bar{r})>(\sinh\psi^{(2)})'(\bar{r})=C_{y_2}(\bar{r}).$$
As $C_u(r)$ is increasing, $C_{u}(r)>C_{y_2}(r)$ for $r\in (\bar{r},a)$. In particular, and 
using  $u'(a)=y_2'(a)$, 
$$0<\int_{\delta}^a(C_{u}(r)-C_{y_2}(r))dr=\int_{\delta}^a\frac{d}{dr}\left(\sinh\psi(r)-
\sinh\psi^{(2)}(r)\right)dr=-\alpha(\delta)$$
in contradiction with (\ref{alfa}).  
\end{proof}

As a consequence of Lemma \ref{lemma-claim} and putting $y_1(a)<u(a)<y_2(a)$, we have:

\begin{theorem} \label{es1} Let $\kappa>0$ and let $u$ be a solution of the  problem  
(\ref{f1})-(\ref{f2}) given by $u=u(r;u_0)$. If $u_0$ is the lowest height of $u$, then
$$u_0-\frac{1}{\kappa u_0}+\sqrt{a^2+\frac{1}{\kappa^2 u_0^2}}<u(a)<u_0+a\frac{\cosh\beta-1}{\sinh\beta}.$$
\end{theorem}

If we compare with (\ref{cotaz22}), the upper bound obtained in Theorem \ref{es1} does not 
depend on $\kappa$ but only on $a$. 
Other source to control the shape of $u$ comes from the integration of $u$. We know 
from (\ref{B11}) and  Lemma \ref{lemma-claim} that 
\begin{equation}\label{y1y2}
\kappa\int_0^a y_1(t)\ dt<\sinh\psi(a)<\kappa\int_0^a y_2(t)\ dt.
\end{equation}
Actually, these inequalities inform us about the volume per unit of lenght that encloses 
each one of the three surfaces together with the support plane $\{x_3=u(a)\}$. The  difference with the 
estimate obtained in Theorem \ref{es1} is that we now obtain a control of the value $\psi(a)=\beta$. For the integrals involving $y_i$, we  write 
$$\int_0^a\left(\sqrt{r^2+m^2}+c\right) dr:=F(c,m)= a c+\frac{a}{2}\sqrt{a^2+m^2}+\frac{m^2}{2}\log(\frac{a+\sqrt{a^2+m^2}}{m}).$$

Then the first inequality in (\ref{y1y2}) yields $\kappa F(u_0-\mu_1,\mu_1)<\sinh\beta$. Thus 
$$\frac{a}{2\kappa u_0}\left(2(\kappa u_0^2-1)+\sqrt{1+a^2\kappa^2 u_0^2}\right)+
\frac{\log(a\kappa u_0+\sqrt{1+a^2\kappa^2 u_0^2})}{2\kappa^2 u_0^2}<\frac{\sinh\beta}{\kappa}.$$
The other inequality in (\ref{y1y2}) says $\sinh\beta<\kappa F(u_0-\mu_2,\mu_2)$. Then
$$\frac{\sinh\beta}{\kappa}<au_0+\frac{a^2\coth\beta}{2}+\frac{a^2\beta}{2\sinh^2\beta}
-\frac{a^2}{\sinh\beta}.$$
In particular, 

\begin{theorem}\label{es2} Fix $\kappa>0$ and let $u=u(r;u_0)$ be a   solution of the  problem  
(\ref{f1})-(\ref{f2}). Then
\begin{equation}\label{laplace}
\frac{\sinh\beta}{a\kappa}+ \frac{a}{\sinh\beta}-\frac{a\coth\beta}{2}-\frac{a\beta}{2\sinh^2\beta}<u_0<\frac{\sinh\beta}{a\kappa}.
\end{equation}
\end{theorem}

\begin{proof} 
The left inequality in (\ref{laplace}) is a consequence of $\sinh\beta<\kappa F(u_0-\mu_2,\mu_2)$; the right inequality comes 
by comparing  the slopes of $y_1$ and $u$ at the point $r=a$: $y_1'(a)<u'(a)$.
\end{proof}

We obtain a new estimate of the solution $u$. 
For this, let us move down the hyperbola $y_2$ until it meets $u$ at $(a,u(a))$. 
We denote by $y_3$ the new position of $y_2$. 

\begin{lemma} The function  $y_3$ satisfies $y_3<u$ on the interval $[0,a)$. 
\end{lemma}
\begin{proof} 
With a similar argument as in Lemma \ref{lemma-claim}, 
 we compare the curvatures of 
$u$ and $y_3$: by (\ref{B2}), we have
$$C_u(a)=\kappa u(a)>\frac{\sin\psi(a)}{a}=C_{y_3}(a).$$
Thus, around the point $r=a$, $y_3<u$. By contradiction, assume that there is $\delta\in (0,a)$ such that 
$y_3(r)<u(r)$ for $r\in(\delta,a)$ and $y_3(\delta)=u(\delta)$. Since $u'(\delta)\geq y_3'(\delta)$, then
$\psi^{(3)}(\delta)\leq \psi(\delta)$. This implies 
\begin{equation}\label{C}
\int_{\delta}^a(C_{y_3}(r)-C_u(r))dr=\sin\psi(\delta)-\sin\psi^{(3)}(\delta)\geq 0.
\end{equation}
Then there would be $\bar{r}\in (\delta,a)$ such that $C_{y_3}(\bar{r})-C_u(\bar{r})>0$.
As $C_u(r)$ is increasing on $r$, $C_u(r)<C_{y_3}(r)$ on $(0,\bar{r})$ and hence also 
throughout $(0,\delta)\subset(0,\bar{r})$. Thus
$$0<\int_0 ^\delta (C_{y_3}(r)-C_u(r))dr=\sin\psi(\delta)-\sin\psi^{(3)}(\delta)\leq 0$$
by (\ref{C}). This contradiction shows the result.
\end{proof}

 As conclusion, we have 
 the estimates:
$$y_3(r)<u(r)\hspace*{1cm}0\leq r<a.$$
$$F(u(a)-a \coth\beta;\mu_2)<\frac{\sinh\beta}{\kappa}.$$
Both inequalities give the next

\begin{theorem} With the same notation as in Theorem \ref{es2}, we have
\begin{equation}\label{meniscus2}
u(a)-a\coth\beta+\sqrt{r^2+\frac{a^2}{\sinh^2\beta}}
<u(r),\hspace*{.5cm}0\leq r<a.
\end{equation}
\begin{equation}\label{ua}
u(a)<\frac{\sinh\beta}{\kappa a}+\frac{a}{2}\coth\beta-\frac{a\beta}{2\sinh^2\beta}.
\end{equation}
\end{theorem}
We point out that the upper bound obtained in (\ref{ua}) does not depend on $u_0$.

\section{Stationary  bands: the case $\kappa<0$}\label{pendent}

This section is devoted to the study of stationary bands  when $\kappa<0$.  We assume in this section that $u_0<0$.

\begin{theorem}\label{t-p}
Let $\Omega_a$ be a strip of the $(x_1,x_2)$-plane, $a>0$.  
Let $u(r;u_0)$ be a solution of the problem (\ref{e1})-(\ref{e2}). Then $u$ is a periodic 
function that vanishes in  an infinite discrete set of points.  The inflections of $u$ are their zeros. Moreover, 
 $u_0\leq u(r)\leq-u_0$, attaining both values  at
exactly the only critical points of $u$.
\end{theorem}

\begin{proof}
From (\ref{B11}), $u'$ is positive near $r=0$ and then, $u$ is strictly increasing on some 
interval $[0,\epsilon)$. As consequence of (\ref{e1}), $u$ is convex around 
$r=0$ and $u$ is convex provided the function $u$  is negative. This implies that $u$ must vanishes at 
some point $r=r_o$. 
From (\ref{cotaz}) the zeroes of $u'$ agree with $u=\pm u_0$ and from 
(\ref{e1}), the inflections agree with the zeroes of $u$. By Theorem \ref{symmetry}, we obtain 
the result.
\end{proof} 

\begin{corollary} Let $S$ be a pendent stationary band. Then 
$S$ is invariant by a group of horizontal translations orthogonal to the rulings. 
\end{corollary}

Using (\ref{cotaz22}) again and Theorem \ref{t-p}, we have

\begin{corollary} Let $\kappa<0$. Then the 
 maximum slope of a solution $u(r;u_0)$ of (\ref{e1})-(\ref{e2}) 
occurs at each  zero 
of $u$ and its value is
$$u'(r_o)=\frac{-u_0}{2-\kappa u_0^2}\sqrt{\kappa^2 u_0^2-4\kappa}.$$
\end{corollary}

We show the existence of pendent stationary bands in the variational 
problem. We need the following

\begin{lemma} \label{pendent-aux} 
Consider $u=u(r;u_0)$ a solution of (\ref{e1})-(\ref{e2}) and 
denote $r_o$ the first zero of $u$. Then 
\begin{equation}\label{pendent-zero}
\sqrt{\frac{-2}{\kappa}}<\sqrt{u_0^2-\frac{2}{\kappa}}<r_o.
\end{equation}
\end{lemma}

\begin{proof}
We consider the hyperbola $y_4$ defined by
$$y_4(r)=\sqrt{r^2+\left(\frac{1}{\kappa u_0}\right)^2}+u_0-\frac{1}{\kappa u_0}.$$
Using the same argument as in  (\ref{B2}), 
the function $u$ is negative in the interval $(0,r_o)$ and then,  $\kappa r u(r)<\sinh\psi(r)<\kappa r u_0$.   
It follows that  $u'(r)<y_4'(r)$. Since $y_4(0)=u(0)$ then $u(r)<y_4(r)$. 
  As $y_4$ meets the $r$-axis at the point
$\sqrt{u_0^2-2/\kappa}$, a comparison between $y_4$ and $u$  gives the desired estimates.
\end{proof}

\begin{theorem}[Existence]\label{pendent-ex}  
Let $\Omega_a$ be a strip 
of the $(x_1,x_2)$-plane, $a>0$.  
Given constants $\kappa<0$ and $\beta$,  there exists a stationary band on 
$\Omega_a$ whose directrix is defined 
by a function $u=u(r;u_0)$, that makes a contact hyperbolic angle 
$\beta$ with the support plane $\{x_3=u(a)\}$. 
\end{theorem}

\begin{proof} 
If $\beta=0$, we take $S=\{x_3=0\}$. Without loss of generality, we now assume 
$\beta>0$. The problem is equivalent to search a solution of (\ref{f1})-(\ref{f2}). For this, we take the initial value problem (\ref{e1})-(\ref{e2}), with $u_0<0$. The 
problem then reduces to find $u_0<0$ such that $u'(a;u_0)=\tanh\beta$. We will search the 
solution in such way that $u$ is negative in its domain. We know by the 
continuity of parameters that $\lim_{u_0\rightarrow 0}u'(a;u_0)=0$. 

On the other hand, we show that $\lim_{u_0\rightarrow -\infty}u'(a;u_0)=1$. For 
this, we know that if $|u_0|$ is sufficiently big, then $a<\sqrt{u_0^2-2/\kappa}<r_o$, $r_o$ 
the first zero of $u(r;u_0)$. It follows from the proof of Lemma \ref{pendent-aux}  that 
$u(r;u_0)<y_4(r)$, for $0<r<a$. Since both functions are negative, we have from 
(\ref{B11}) that
$$\sinh\psi(a)=\kappa\int_0^{a} u(t) dt>\kappa\int_0^a y_4(t)\ dt=F(u_0-\frac{1}{\kappa u_0},\frac{1}{\kappa u_0}) 
\longrightarrow +\infty,$$
when $u_0\rightarrow-\infty$. Thus $u'(a;u_0)=\tanh\psi(a)\rightarrow 1$, as $u_0\rightarrow-\infty$. By the continuity of parameters, we conclude the existence of a number $u_0<0$ with the desired condition of Theorem \ref{pendent-ex}.
\end{proof}

Finally, we establish
 some estimates of the solution $u$ for the problem (\ref{f1})-(\ref{f2}). 
By the periodicity of $u$, we restrict to the interval $[0,r_o]$ and that $a\leq r_o$ (for example, this condition holds if 
$a<\sqrt{-2/\kappa}$, see Lemma \ref{pendent-aux}). We use the function  
$y_4$. In particular, we know for $0<r\leq a$
$$u(r)<y_4(r),\hspace*{1cm}
\int_0^a u(t)\ dt<\int_0^a y_4(t)\ dt.$$
A somewhat similar argument as in Theorems \ref{es1} and \ref{es2}, we 
conclude, respectively,
$$u(a)<u_0-\frac{1}{\kappa u_0}+\sqrt{a^2+\frac{1}{\kappa^2 u_0^2}}$$
$$\frac{\sinh\beta}{\kappa}<\frac{a}{2\kappa u_0}\left(2(\kappa u_0^2-1)+\sqrt{1+a^2\kappa^2 u_0^2}\right)+
\frac{\log(a\kappa u_0+\sqrt{1+a^2\kappa^2 u_0^2})}{2\kappa^2 u_0^2}.$$

 For pendent stationary bands we have not a result of monotonicity with respect to the parameters. This is due to the periodicity of solutions. At this state, we  can only 
assure the following results of monotonicity on a certain interval around $r=0$: 
\begin{enumerate}
\item If $\kappa_1<\kappa_2$, then $u(r;u_0,\kappa_1)>u(r;u_0,\kappa_2)$ (Theorem \ref{mono-k}).
\item Let $\delta>0$. Then $u(r;u_0-\delta)+\delta\geq u(r;u_0)$ (Theorem \ref{delta}).
\end{enumerate}


\end{document}